\documentclass[12pt]{article}
\usepackage{amsmath,amssymb,amsthm, amsfonts, bm}
\usepackage[mathscr]{eucal}
\usepackage{stackengine}
\usepackage{accents}
\usepackage{rotate,graphics,epsfig}
\usepackage{color}

\newtheorem{The}{Theorem}
\newtheorem{Exa}[The]{Example}
\newtheorem{Cor}[The]{Corollary}

\newtheorem{Pro}[The]{Proposition}

\theoremstyle{definition}
\newtheorem{Def}[The]{Definition}
\newtheorem{Rem}[The]{Remark}
\numberwithin{equation}{section}
\numberwithin{The}{section}
\numberwithin{figure}{section}

\usepackage[utf8]{inputenc}
\usepackage[T1]{fontenc}

\newcommand{\be}{\begin{eqnarray}}
\newcommand{\ee}{\end{eqnarray}}

\newcommand{\by}{\begin{eqnarray*}}
\newcommand{\ey}{\end{eqnarray*}}

\usepackage{color}

\begin{document}
\title{Toward a Hazard Rate Framework for Regular and Rapid Variation}
\author{
Haijun Li
\footnote{{\small\texttt{lih@math.wsu.edu}}, Department of Mathematics and Statistics, Washington State University, 
	Pullman, WA 99164, U.S.A.}
}
\date{April 2025}
\maketitle

\begin{abstract}

Regular and rapid variation have been extensively studied in the literature and applied across various fields, particularly in extreme value theory. In this paper, we examine regular and rapid variation through the lens of generalized hazard rates, with a focus on the behavior of survival and density functions of random variables. Motivated by the von Mises condition, our hazard rate based framework offers a unified approach that spans from slow to rapid variation, providing in particular new insights into the relationship between hazard rate functions and the right tail decays of random variables.	
	
\end{abstract}

\section{Introduction}
\label{Instroduction}

Let $X_1, \dots, X_n$ be $n$ independent and identically distributed (i.i.d.) real-valued random samples with distribution function $F$, or equivalently,  survival function $\overline{F}=1-F$, having the right support point at $+\infty$. The univariate extreme-value theory concerns the limiting behaviors of the sample extremes,  under appropriate normalization, as $n$ goes to infinity. The well-known Fisher–Tippett–Gnedenko theorem states that the properly affine-normalized maxima converge in distribution to a generalized extreme-value distribution \cite{Haan-Ferreira06}, which essentially consists of two main types of distributions: the Fréchet and the Gumbel distributions. The Gnedenko–de Haan theorem (see \cite{Gnedenko43, de70, Haan-Ferreira06}) further refines this asymptotic result by characterizing the domains of attraction of these limit laws in terms of the underlying sampling distribution function $F$ as follows.
\begin{enumerate}
	\item The limiting distribution is the Fr\'echet distribution $H_+(x; \alpha) = \exp\{-x^{-\alpha}\}$, $x\ge 0$, $\alpha >0$, if and only if the sampling survival function $\overline{F}$ is {\em regularly varying} with tail parameter $-\alpha<0$, i.e.,
	\begin{equation}\label{eq-1}
	\lim_{t\to \infty}\frac{\overline{F}(tx)}{\overline{F}(t)}= x^{-\alpha}, \ x>0.
	\end{equation}
	\item The limiting distribution is the Gumbel distribution $H_0(x) = \exp\{-e^{-x}\}$, $x\in \mathbb{R}$, if and only if the sampling survival function $\overline{F}$ is {\em rapidly varying} in the sense that
		\begin{equation}\label{eq-2}
	\lim_{t\to \infty}\frac{\overline{F}(t+xR(t))}{\overline{F}(t)}=\exp\{-x\}, x\in \mathbb{R},
		\end{equation}
	where $R(t):=\int_{t}^{\infty}(1-F(x))dx/(1-F(t))$ is known as the {\em reciprocal hazard rate} of the function $\int_t^\infty\overline{F}(x)dx$. 
\end{enumerate}
That is, for the maximum of i.i.d. random samples to converge to an extreme-value limiting distribution under affine normalization, the sampling distribution must be either regularly varying \eqref{eq-1} or rapidly varying \eqref{eq-2}.

Univariate regular and rapid variation have been extensively studied in the literature and have found diverse applications in fields such as differential equations and number theory, among others; see \cite{BGT1987, Omey1981, Omey1997} and references therein. In this paper, we examine univariate regular and rapid variation through the lens of hazard rates, presenting the relationships between right tail decays and asymptotic hazard rates for regular variation and rapid variation.

Let $f: \mathbb{R}\to \mathbb{R}_+$ be non-negative and differentiable. The generalized hazard rate is defined as $h(t) = |f'(t)|/f(t)$, $t\in \mathbb{R}$. The von Mises condition \cite{Mises36} provides a sufficient condition that implies regular variation. Let the sampling distribution $F$ have a positive derivative on $[x_0, \infty)$, and $h_F(x) = F'(x)/(1-F(x))$ denote the hazard rate of $F$. 
		\begin{equation}
	\mbox{If}\ \lim_{t\to \infty}t\,h_F(t)=\alpha \ge 0,\ 
\mbox{then  $\overline{F}$ is regularly varying with tail parameter $-\alpha$.}		
	\label{intro-eq-5}
\end{equation}
The von Mises condition is also necessary for regular variation, provided that the density of $F$ is ultimately decreasing \cite{Landau, Geluk1987, Resnick07}. Note, however, that the von Mises condition does not adequately account for rapid variation. A well-known counterexample is the standard normal distribution $\Phi(x)$, which is rapidly varying, yet satisfies $ \lim_{t\to \infty}t\,h_\Phi(t) =\infty$, at the rate of $t^2$. The issue arises from the use of the scaling function 
$t$ in the von Mises condition \eqref{intro-eq-5}, which accelerates the limiting process and leads to a divergence to infinity, thereby obscuring information about the rate of convergence. In this paper, we address this by replacing the scaling function 
$t$ with 
$R(t)$, as defined for rapid variation in \eqref{eq-2}, to derive a complete analog of the von Mises-type condition suitable for the rapidly varying case.

The methodological similarity between univariate regular and rapid variation plays a key role in our study of multivariate extremes, which is based on the tail dependence of copulas (\cite{JLN10, NJL12, HJL12, LH14, LW2013}). Copula-based tail dependence captures the scale-free extremal dependence structure of multivariate distributions, independent of the univariate margins, which may be either regularly or rapidly varying. When the univariate margins are regularly varying, tail dependence is shown to be equivalent to multivariate regular variation (\cite{Li2009, LS2009, LH13}). Similarly, when the univariate margins are rapidly varying, tail dependence can be used to study multivariate rapid variation (\cite{JL19, Li2021}).

The paper is organized as follows. Section \ref{univariate RV section} discusses regular variation, with a focus on Karamata's representation and Karamata's theorem, which share conceptual similarities with the von Mises condition. Section \ref{univeriate rapid variation} is devoted to rapid variation, drawing extensively on results from \cite{Omey2013}. Section 4 concludes the paper. Throughout, the terms "increasing" and "decreasing" are used in the weak sense, and measurability is assumed without explicit mention unless otherwise stated. Two non-zero measurable functions $f, g: \mathbb{R}\to \mathbb{R}$ are said to be tail equivalent at $\infty$, denoted by $f(t)\sim g(t)$, if $f(t)/g(t)\to 1$ as $t\to \infty$.

\section{Regular Variation}
\label{univariate RV section}

The maximum domain of attraction for a Fr\'echet distribution consists of regularly varying distributions. Regular variation, initiated by Jovan Karamata in 1930, deals with the following limiting process for a function $f:\mathbb{R}_+\to \mathbb{R}_+$:
\begin{equation}
	\frac{f(tx)}{f(t)}\to g(x), \ t\to \infty, \ \forall\ x>0. 
\label{RV-eq-1} 
\end{equation}
The limiting function $g(\cdot)$, if exists, is clearly multiplicative, and such a group invariance property yields via Cauchy's lemma (or its multiplicative version, Hamel's lemma) the explicit expression $g(x) = x^\rho$, $\rho\in \mathbb{R}$.

\begin{Def}
	\label{RV-d-1}
A measurable function $f:\mathbb{R}_+\to \mathbb{R}_+$ is said to be regularly varying at $\infty$ with tail parameter $\rho\in \mathbb{R}$, denoted by $f\in \mbox{RV}_\rho$, if 
\begin{equation}
\lim_{t\to \infty}\frac{f(tx)}{f(t)}=x^\rho,\ \forall\ x>0. 
	\label{RV-eq-2} 
\end{equation}
In particular, $f$ is said to be slowly varying at $\infty$ if $f\in \mbox{RV}_0$. 
\end{Def}

\begin{Rem}
	\label{RV-r-1}
	\begin{enumerate}
		\item In addition to Lebesgue-measurability, there are several other conditions, such as the Baire property, under which \eqref{RV-eq-1} and \eqref{RV-eq-2} are equivalent \cite{BGT1987}. In general, however, assuming the Axiom of Choice, there exist infinitely many non-power functions $g$, e.g., Hamel functions, that satisfy \eqref{RV-eq-1} in this type of Cauchy-Hamel problems. The detailed proof on the equivalence of \eqref{RV-eq-1} and \eqref{RV-eq-2}, based on a Steinhaus theorem, can be found in \cite{BGT1987}, that provides a comprehensive and in-depth study on univariate regular variation, upon which this section is based. 
		\item One can define regular variation toward any point in $\overline{\mathbb{R}}$. For example, a function $f: \mathbb{R}_+\to \mathbb{R}_+$ is said to be regularly varying at $0$ with tail parameter $\rho$ if $f(x^{-1})\in \mbox{RV}_{-\rho}$.
		\item If $f\in \mbox{RV}_\rho$, then $\ell(x):= f(x)x^{-\rho}\in \mbox{RV}_0$ is slowly varying. That is, for any $f\in \mbox{RV}_\rho$, one can always write $f(x) =x^{\rho}\ell(x)$ where $\ell \in \mbox{RV}_0$. 				
	\end{enumerate}
\end{Rem}

A Fr\'echet distribution $H_+(x; \gamma) = \exp\{-x^{-1/\gamma}\}$, $x\ge 0$, $\gamma >0$, has the survival function $1-H_+(x; \gamma)\in \mbox{RV}_{-1/\gamma}$. Moreover, the maximum of i.i.d. random samples converges to the Fr\'echet distribution $H_+(x; \gamma)$ under affine normalization with $\gamma>0$ if and only if the survival function $1-F\in \mbox{RV}_{-1/\gamma}$, excluding the distributions with slowly varying survival functions as sampling distributions for the maximum domains of attraction under affine-normalization. Note, however, that many fundamental properties of regular variation stem from underlying slow variation, which is interpreted as having the zero derivative at infinity.

\begin{Pro}
	\label{RV-pro-1}
	Let $h(x) := \log\ell(e^x)$, $x\in \mathbb{R}$. Then $\ell\in \mbox{RV}_0$ if and only if
	\begin{equation}
		\label{RV-eq-3}
		\lim_{x\to \infty}(h(x+u)-h(x))=0, \ \forall\ u\in \mathbb{R}. 
	\end{equation}
\end{Pro}

\noindent
{\sl Proof.}
The equivalence is immediate, due to the fact that $h$ is additive if and only if $\ell$ is multiplicative. 
\hfill $\Box$

Suppose that $\ell(x)$ is twice differentiable. It then follows from Proposition \ref{RV-pro-1} that the derivative at infinity, $h'(\infty) = \lim_{x\to \infty}h'(x) =0$, if and only if $\ell\in \mbox{RV}_0$. Since this observation is facilitated by interchanging iterated limiting processes for $x$ and differentiation, the local uniform convergence of \eqref{RV-eq-3} is implicitly assumed for slow variation, in the essence of the Moore-Osgood theorem. In fact, the local uniform convergence holds for regular variation.

\begin{The}
	\label{RV-the-1}
	(Karamata's Uniform Convergence Theorem \cite{BGT1987, Resnick07}) 
	\begin{enumerate}
		\item If $\ell\in \mbox{RV}_0$, then 
		\[\frac{\ell(tx)}{\ell(t)}\to 1
		\]
		locally uniformly in $x$ on $(0,\infty)$. 
		\item If $f\in \mbox{RV}_\rho$, $\rho\in \mathbb{R}$, then 
		\[\frac{f(tx)}{f(t)}\to x^{\rho}
		\]
		locally uniformly in $x$ on $(0,\infty)$. For $\rho<0$, the uniform convergence also holds on $(a, \infty]$, $a>0$.
	\end{enumerate}
\end{The}

\noindent
{\sl Proof.}
Several proofs are detailed in \cite{BGT1987}, and some of them provide the direct proofs for the local uniform convergence of \eqref{RV-eq-3}. For $\rho<0$ in (2), one needs to make use of the Alexandroff uncompactification so that $(a,\infty]$ becomes relative compact. 
\hfill $\Box$

\begin{Rem}
	\label{RV-r-2} 
	\begin{enumerate}
		\item 	The uniform convergence of regular variation may fail without imposing some regularity condition, such as measurability, and see \cite{BGT1987} for a counterexample that is based on a Hamel basis, whose existence relies on Zorn's lemma (equivalent to the Axiom of Choice).
		\item The study on regular variation, in its basic form of first-order convergences, is known as Karamata Theory. The convergence rates of regular variation, part of de Haan Theory, are also extensively studied (see \cite{BGT1987} for detail). In particular, one may also study very slow variation, a subclass of slow variation, in which the convergence in \eqref{RV-eq-3} is replaced by the convergence at a certain rate
			\begin{equation}
			\label{RV-eq-35}
		\lim_{x\to \infty}\phi(x)(h(x+u)-h(x))=0, \ \forall\ u\in \mathbb{R},
		\end{equation}
		where the rate function $\phi:\mathbb{R}_+\to (0,\infty)$ is increasing. Under the measurability assumption, \eqref{RV-eq-35} also holds locally uniformly in $u$. 
	\end{enumerate}
\end{Rem}

It follows immediately from Theorem \ref{RV-the-1} that a regularly varying function has monotone tail equivalents.

\begin{Pro}
	\label{RV-pro-2} Suppose that  $f: \mathbb{R}_+\to \mathbb{R}_+$ is locally bounded, and  $f\in \mbox{RV}_\rho$, $\rho\in \mathbb{R}$.
	\begin{enumerate}
		\item If $\rho>0$, then the increasing regularly varying function
		\[\bar{f}(t) :=\sup\{f(x): 0\le x\le t\}\sim f(t),\ t\to \infty. 
		\]
		\item If $\rho<0$, then the decreasing regularly varying function 
			\[ \underaccent{\bar}{{f}}(t):=\inf\{f(x): 0\le x\le t\}\sim f(t),\ t\to \infty. 
		\]
	\end{enumerate}
\end{Pro}

\noindent
{\sl Proof.} Consider the case that $\rho>0$. 
Observe that $\bar{{f}}(t)=\sup_{0\le x\le 1}f(tx)$, $t>0$. By Theorem \ref{RV-the-1} (2),
\[\frac{\bar{{f}}(t)}{f(t)} = \sup_{0\le x\le 1}\frac{f(tx)}{f(t)}\to \sup_{0\le x\le 1}x^{\rho} = 1; 
\]
that is, $\bar{{f}}(t)\sim f(t)$. The case that $\rho<0$ follows immediately by using $1/f(t)$. 
\hfill $\Box$

Local uniform convergence yields various representations for slow and regular variations.

\begin{The}
	\label{RV-the-2}
	(Karamata's Representation \cite{BGT1987, Resnick07}) 
	\begin{enumerate}
		\item A function $\ell\in \mbox{RV}_0$ if and only if
		\[\ell(x) = c(x)\exp\left\{\int_{x_0}^x\epsilon(t)dt/t\right\}
		\]
		for some constant $x_0>0$, where the functions $c, \epsilon: \mathbb{R}_+\to \mathbb{R}$ have the limits $\lim_{x\to \infty}c(x) = c\in (0,\infty)$ and $\lim_{x\to \infty}\epsilon(x) = 0$. 
		\item In general, a function $f\in \mbox{RV}_\rho$, $\rho\in \mathbb{R}$, if and only if
		\[f(x) = c(x)\exp\left\{\int_{x_0}^x\rho(t)dt/t\right\}
		\]
		for some constant $x_0>0$, where the function $c, \rho: \mathbb{R}_+\to \mathbb{R}$ have the limits $\lim_{x\to \infty}c(x) = c\in (0,\infty)$ and $\lim_{x\to \infty}\rho(x) = \rho$. 
	\end{enumerate}
\end{The}

\noindent
{\sl Proof.}
The representation of slow variation is shown to be equivalent to the uniform convergence of slow variation (see Theorem \ref{RV-the-1} (1)), and the proof can be found in \cite{BGT1987}. The representation for regular variation (2) follows immediately from (1) and the fact that $f(x) = x^{\rho}\ell(x)$, $x\in \mathbb{R}_+$. 
\hfill $\Box$

\begin{Rem}
	\label{RV-r-3} 
	\begin{enumerate}
		\item The functions in Karamata's representation are not uniquely determined, and the functions $c(\cdot)$, $\epsilon(\cdot)$ and $\rho(\cdot)$ may be adjusted on any compact sets accordingly.  
		The function $c(\cdot)$ in Karamata's representation can be taken as being ultimately bounded. A slowly varying function $\ell(\cdot)$ can be taken as being ultimately locally bounded.
	\item The function $\epsilon(\cdot)$ may be assumed to be arbitrarily smooth, and as a result, the function
	\begin{equation}
		\label{RV-eq-4}
		\ell_1(x) = c\exp\left\{\int_{x_0}^x\epsilon(t)dt/t\right\},\ x\ge x_0, 
		\end{equation}
	known as the normalized slowly varying function of $\ell(\cdot)$, is also arbitrarily smooth. Note that $\ell(x)\sim \ell_1(x)$ as $x\to \infty$. 
	\item If $\ell(\cdot)$ is ultimately increasing (ultimately decreasing), then $\epsilon(\cdot)$ may be taken as being ultimately non-negative (ultimately non-positive). 
	\item Since the limit at infinity is of interest for regular variation, the uniform convergence theorem and Karamata's representation hold for the entire class with respect to tail equivalence. That is, if $f(x)\sim g(x)$ as $x\to \infty$ and any one of Theorems \ref{RV-the-1} and \ref{RV-the-2} holds for $f$, it also holds for $g$. 
	\end{enumerate}
\end{Rem}

Typical slowly varying functions involve multiplicative factors such as powers of iterated natural logarithms and $\exp\big\{a(\log x)^b\big\}$, $a\in \mathbb{R}$, $0<b<1$. 
Various specific examples of slow and regular variations can be found in \cite{BGT1987, Haan-Ferreira06, Resnick07}, which also include in detail the properties of regular variation. These properties can be derived from Definition \ref{RV-d-1} and Karamata's representation, Theorem \ref{RV-the-2}. Some useful properties  are summarized below.

\begin{Pro}
	\label{RV-pro-25}
\begin{enumerate}
	\item If $f\in \mbox{RV}_\rho$, $\rho\in \mathbb{R}$, then $f^\alpha\in \mbox{RV}_{\rho \alpha}$, $\forall\ \alpha\in \mathbb{R}$. 
	\item If $f\in \mbox{RV}_\rho$, $\rho\in \mathbb{R}$, then $\lim_{x\to \infty}\log f(t)/\log t = \rho$, implying that
	\[	\lim_{x\to \infty}f(t) = \left\{\begin{array}{ll}
		0  & \mbox{if $\rho< 0$}\\
		\infty  & \mbox{if $\rho>0$.}
	\end{array}
	\right.
	\]
	\item If $f\in \mbox{RV}_\rho$, $\rho\in \mathbb{R}$, then for any $\epsilon>0$, there exists $t_0=t_0(\epsilon)$, such that for any $x\ge 1$ and $t\ge t_0$, 
	\begin{equation}
		\label{RV-eq-5}
		(1-\epsilon)x^{\rho-\epsilon}<\frac{f(tx)}{f(t)}< (1+\epsilon)x^{\rho+\epsilon}.
	\end{equation}
\item If $f\in \mbox{RV}_\rho$, $\rho> 0$, and $f$ is increasing and $\lim_{x\to \infty}f(x) = \infty$, then the left-continuous inverse $f^{\leftarrow}(s)=\inf\{t:f(t)\ge s\}\in \mbox{RV}_{1/\rho}$. If $f\in \mbox{RV}_0$, then left-continuous inverse $f^{\leftarrow}$ is {\em rapidly varying} in the sense of de Haan:
\begin{equation}
	\label{RV-eq-6}
	\lim_{t\to \infty}\frac{f^{\leftarrow}(tx)}{f^{\leftarrow}(t)} = \left\{\begin{array}{lll}
		0  & \mbox{if $x < 1$}\\
		1  & \mbox{if $x=1$}\\
		\infty & \mbox{if $x>1$.}
	\end{array}
	\right.
\end{equation}
	\item If $f_1\in \mbox{RV}_{\rho_1}$ and $f_2\in \mbox{RV}_{\rho_2}$, $\rho_1, \rho_2\in \mathbb{R}$, than $f_1f_2\in \mbox{RV}_{\rho_1+\rho_2}$, $f_1+f_2\in \mbox{RV}_{\max\{\rho_1+\rho_2\}}$. If moreover $\lim_{x\to \infty}f_2(x) = \infty$, then the composition $f_1\circ f_2\in \mbox{RV}_{\rho_1\rho_2}$. 
	\item If two increasing functions $f_1\in \mbox{RV}_{\rho}$ and $f_2\in \mbox{RV}_{\rho}$, $\rho>0$, then for any constant $c>0$, 
	\[f_1(t)\sim cf_2(t), \ t\to \infty,\ \mbox{if and only if},\ f_1^{\leftarrow}(x)\sim c^{-1/\rho} f_2^{\leftarrow}(x),\ x\to \infty. 
	\]
\end{enumerate}
\end{Pro}

Proposition \ref{RV-pro-25} (4) and (6) also hold for the case of decreasing regularly varying functions. Define the right-continuous inverse $g^{\rightarrow}(x):=\inf\{y: g(y)\le x\}$ for a decreasing function $g: \mathbb{R}_+\to \mathbb{R}_+$. Note that $g^{\rightarrow}(x) \le t$ if and only if $x \ge g(t)$, and $g^{\rightarrow}(x)$ is decreasing. 
\begin{Cor}
	\label{RV-pro-3} 
	\begin{enumerate}
		\item If $f\in \mbox{RV}_\rho$, $\rho< 0$, and $f$ is decreasing and $\lim_{x\to \infty}f(x) = 0$, then the right-continuous inverse $f^{\rightarrow}(t^{-1})\in \mbox{RV}_{-1/\rho}$.
		\item If two decreasing functions $f_1\in \mbox{RV}_{\rho}$ and $f_2\in \mbox{RV}_{\rho}$, $\rho<0$, then for any constant $c>0$, 
		\[f_1(t)\sim cf_2(t), \ t\to \infty,\ \mbox{if and only if},\ f_1^{\rightarrow}(x^{-1})\sim c^{1/\rho} f_2^{\rightarrow}(x^{-1}),\ x\to \infty. 
		\]
	\end{enumerate}
\end{Cor}

\noindent
{\sl Proof.} (1)
Let $\hat{f}(t) = 1/f(t)$, $t>0$. Observe that $\hat{f}\in \mbox{RV}_{-\rho}$, $-\rho>0$,  and $\hat{f}$ is increasing and $\lim_{t\to \infty}\hat{f}(t)=\infty$. By Proposition \ref{RV-pro-25} (4), $\hat{f}^{\leftarrow}\in \mbox{RV}_{-1/\rho}$. Since
\[\hat{f}^{\leftarrow}(t^{-1}) = \inf\big\{s: \hat{f}(s)\ge t^{-1}\big\}=\inf\big\{s: f(s)\le t\big\}=f^{\rightarrow}(t),
\]
one has that $f^{\rightarrow}(x^{-1})\in \mbox{RV}_{-1/\rho}$.

(2) Similarly, let $\hat{f}_i(t) = 1/f_i(t)$, $t>0$, $i=1,2$. It follows from Proposition \ref{RV-pro-25} (6) that $f_1(t)\sim cf_2(t)$ is equivalent to that $\hat{f}_1^{\leftarrow}(x)\sim c^{1/\rho} \hat{f}_2^{\leftarrow}(x)$, which is, in turns, equivalent to $f_1^{\rightarrow}(x^{-1})\sim c^{1/\rho} f_2^{\rightarrow}(x^{-1}),\ x\to \infty$. 
\hfill $\Box$

It is worth mentioning that Proposition \ref{RV-pro-2} (4) and (6) hold for ultimate increasing functions, and Corollary \ref{RV-pro-3} holds for ultimately decreasing functions. 
Theorems \ref{RV-the-1} and \ref{RV-the-2} yield Karamata's theorem on integration of regularly varying functions.

\begin{The}
	\label{RV-the-3}
	(Karamata's Theorem \cite{BGT1987, Haan-Ferreira06, Resnick07}) 
	Let $f: \mathbb{R}_+\to \mathbb{R}_+$ be ultimately locally bounded.
	\begin{enumerate}
		\item If $f\in \mbox{RV}_\rho$ with $\rho\ge -1$, then $\int_{t_0}^tf(x)dx\in \mbox{RV}_{\rho+1}$ and 
		\begin{equation}
			\label{RV-eq-7}
			\lim_{t\to \infty}\frac{tf(t)}{\int_{t_0}^tf(x)dx} = \rho+1.
		\end{equation}
	If $f\in \mbox{RV}_\rho$ with $\rho< -1$ or $f\in \mbox{RV}_{\rho}$ with $\rho=-1$ and $\int_t^\infty f(x)dx < \infty$, then $\int_t^\infty f(x)dx\in \mbox{RV}_{\rho+1}$ and 
		\begin{equation}
		\label{RV-eq-8}
		\lim_{t\to \infty}\frac{tf(t)}{\int_t^\infty f(x)dx} = -\rho-1.
	\end{equation}
\item Conversely, if $f$ satisfies \eqref{RV-eq-7} with $-1<\rho<\infty$, then $f\in \mbox{RV}_\rho$; if $f$ satisfies \eqref{RV-eq-8} with $-\infty<\rho<-1$, then $f\in \mbox{RV}_\rho$.
	\end{enumerate}
\end{The}

\begin{Rem}
		\label{RV-r-4} 
If $f(\cdot)$ is the density of a positive random variable $X$, then $f(t)/\int_t^\infty f(x)dx$ is the hazard rate, and $f(t)/\int_0^tf(x)dx$ is the reverse hazard rate of $X$. With such interpretations, both \eqref{RV-eq-7} and \eqref{RV-eq-8} provide von Mises-type conditions for regularly varying densities and regularly varying survival functions. 
\end{Rem}

\begin{Exa}
	\label{RV-exa-4}\rm
	Karamata's theorem can be used to strengthen Proposition \ref{RV-pro-2}. Consider $g\in \mbox{RV}_\rho$, $\rho>0$. Define
	\[\bar{g}(t) = \int_1^tg(s)ds/s, \ t\ge 1,
	\]
	which is absolutely continuous and ultimately strictly increasing. 
	Since $t^{-1}g(t)\in \mbox{RV}_{\rho-1}$, it follows from Karamata's theorem that $g(t)/\bar{g}(t)\to \rho$, and hence $g(t)\sim \rho \bar{g}(t)$ as $t\to \infty$. That is, one can find a tail equivalent of $g$, that is  ultimately strictly increasing, and  absolutely continuous with derivative $\rho\bar{g}'(t)\in \mbox{RV}_{\rho-1}$.
	
	Similarly, if $g\in \mbox{RV}_\rho$, $\rho<0$, then $g$ is tail equivalent to a strictly decreasing differentiable function with a negative derivative that is regularly varying with tail parameter $\rho-1$. 
 \hfill $\Box$
\end{Exa}

\begin{Exa}
	\label{RV-exa-3}\rm
Consider a density function $f(x) = x^{-1}(\log x)^{-2}\in \mbox{RV}_{-1}$, $x>e$. The corresponding survival function $\overline{F}(t) = 1/\log x$ is slowly varying and \eqref{RV-eq-8} is satisfied. \hfill $\Box$
\end{Exa}
 
If $f\in \mbox{RV}_{-1}$, then both $\int_{t_0}^tf(x)dx$ and $\int_t^\infty f(x)dx$ are slowly varying and satisfy \eqref{RV-eq-7} and \eqref{RV-eq-8}, respectively. Note that 	
when $\rho=-1$, \eqref{RV-eq-7} and \eqref{RV-eq-8} are only two necessary conditions for $f\in \mbox{RV}_{-1}$, due to the zero limit that provides no information on asymptotic rates. 
It is worth mentioning, however, that when $\rho = -1$, \eqref{RV-eq-7} still implies that $\int_{t_0}^tf(x)dx\in \mbox{RV}_0$, and \eqref{RV-eq-8}  implies that $\int_{t}^\infty f(x)dx\in \mbox{RV}_0$, via the von Mises condition \cite{Resnick07}.

\begin{The}
	\label{RV-the-5} (von Mises \cite{Mises36}) Let $g:\mathbb{R}_+\to \mathbb{R}_+$ be differentiable. If the von Mises condition
		\begin{equation}
		\label{RV-eq-20}
		\lim_{t\to \infty}t\,
		\frac{g'(t)}{g(t)} = \rho\in \mathbb{R}, 
	\end{equation}
holds, then $g\in \mbox{RV}_\rho$. 
\end{The}

Conversely, the monotonicity is needed. 

\begin{The}
\label{RV-the-4} (Landau \cite{Landau}) If a differentiable function $g\in \mbox{RV}_{\rho}$ and $g'$ is ultimately monotone, then the von Mises condition \eqref{RV-eq-20} holds, and in addition, 
$g'\in \mbox{RV}_{\rho-1}$ if $\rho>0$ and $-g'\in \mbox{RV}_{\rho-1}$ if $\rho<0$.
\end{The}

\begin{Rem}
	\label{RV-r-5} 
	\begin{enumerate}
		\item Let $g(t)=\mathbb{P}(X>t)$ denote the survival function  of a positive random variable $X$, and then $-g'$ is its density. If $g\in \mbox{RV}_{\rho}$, $\rho<0$, and the density is ultimately decreasing, then the von Mises condition \eqref{RV-eq-20} holds for the survival function of $X$, and its density is regularly varying. 
		\item Theorems \ref{RV-the-5} and \ref{RV-the-4} can be summarized,  in terms of hazard rates, as follows. If $g:\mathbb{R}\to \mathbb{R}$ is ultimately positive with ultimately monotone derivative $g'(t)$, then 
			\begin{equation}
			\lim_{t\to \infty}t\,\frac{g'(t)}{g(t)}=	\left\{\begin{array}{ll}
				0 & \mbox{if and only if $g$ is slowly varying at $\infty$;}\\
				\rho\in \mathbb{R}\backslash \{0\}  & \mbox{if and only if $g$ is regularly varying at $\infty$.}
			\end{array}
			\right.
			\label{RV-eq-15}
		\end{equation}
	\end{enumerate}
\end{Rem}

\section{Rapid Variation}
\label{univeriate rapid variation}

A regularly varying function $f\in \mbox{RV}_\rho$, $\rho\in \mathbb{R}$, behaves asymptotically like a power function with power $\rho$, and thus, if $\rho\ne 0$, then $h(t):=f(e^{t})$ is rapidly varying in the sense of de Haan:
\begin{equation}
\lim_{t\to \infty}\frac{h(tx)}{h(t)}=	\left\{\begin{array}{ll}
		\infty & \mbox{if $0<x<1$}\\
		0  & \mbox{if $x>1$,}
	\end{array}
  \right.
	\label{Rapid-eq-1}
\end{equation}
for $\rho<0$; and for $\rho>0$, 
\begin{equation}
	\lim_{t\to \infty}\frac{h(tx)}{h(t)}=	\left\{\begin{array}{ll}
		0 & \mbox{if $0<x<1$}\\
		\infty  & \mbox{if $x>1$.}
	\end{array}
	\right.
	\label{Rapid-eq-2}
\end{equation}
The class  of rapidly varying functions satisfying \eqref{Rapid-eq-1} (or \eqref{Rapid-eq-2}) is denoted by   $\mbox{R}_{-\infty}$ (or $\mbox{R}_\infty$). 
To analyze the limiting behaviors, decrease the process rate for the numerator and reformulate the limiting ratio for a regularly varying function $f\in \mbox{RV}_\rho$, $\rho\in \mathbb{R}$, thereby obtaining an exponential-type limit,
\begin{equation}
	\label{Rapid-eq-3}
	\lim_{t\to \infty}\frac{h(t+x)}{h(t)}= \lim_{t\to \infty}\frac{f(e^{t+x})}{f(e^t)} = \exp\{\rho x\},\ \forall\ x\in \mathbb{R}. 
\end{equation}
That is, $h:\mathbb{R}\to \mathbb{R}_+$ is known as a function with a long right tail, denoted by $h\in \mbox{L}_\rho$. Therefore, $f\in \mbox{RV}_\rho$, if and only if $f(e^{t}) \in \mbox{L}_\rho$. A general class of functions that cover both \eqref{Rapid-eq-3} and functions with exponential growths or decays in $\mbox{R}_{\pm \infty}$ is defined as follows.

\begin{Def}\label{Rapid-def-1}
	An ultimately positive, measurable function $f:\mathbb{R}\to \mathbb{R}$ is said to be in $\Gamma_\alpha(g)$, denoted by $f\in \Gamma_\alpha(g)$, if for $\alpha\in \mathbb{R}$, there exists a measurable and positive function $g$, such that
	\begin{equation}
		\label{Rapid-eq-4}
		\lim_{t\to \infty}\frac{f(t+xg(t))}{f(t)} = \exp\{\alpha x\},\ \forall\ x\in \mathbb{R}. 
	\end{equation}
\end{Def}

\begin{Rem}\label{Rapid-rem-1}
	\begin{enumerate}
		\item If $\alpha > 0$ ($\alpha < 0$), then $f\in \Gamma_\alpha(g)$ if and only if $f\in \Gamma_1(g/\alpha)$ $\big(f\in \Gamma_{-1}(g/(-\alpha))\big)$. The detailed discussion on $\Gamma_\alpha(g)$ and convergence rates can be found in \cite{Omey2013}, upon which this section is based.
		\item If $g(t) = 1$, then $\Gamma_\alpha(g)=\mbox{L}_{\alpha}$; that is, $f\in \Gamma_\alpha(1)$ if and only if $f(\log t)\in \mbox{RV}_{\alpha}$. 
		 In addition, it follows from the equivalence of \eqref{RV-eq-1} and \eqref{RV-eq-2} that if $g(t)/t\to \theta> 0$, then the exponential-type limit \eqref{Rapid-eq-4} does not hold for any finite $\alpha$. That is, in order for \eqref{Rapid-eq-4} to hold for a finite $\alpha$, $g(t)/t\to 0$. 
	\end{enumerate}
\end{Rem}

It is known from \cite{Geluk1987, BGT1987} that if \eqref{Rapid-eq-4} holds for a monotone function $f\in \Gamma_{\alpha}(g)$, then the corresponding function $g$ is Beurling varying; that is, 
\begin{equation}
	\label{Rapid-eq-5}
	g(t)/t\to 0, \ g\in \Gamma_0(g),
\end{equation}
and moreover, 	
\begin{equation}
	\label{Rapid-eq-6}
	\lim_{t\to \infty}\frac{g(t+xg(t))}{g(t)} = 1,\ \forall\ x\in \mathbb{R},
\end{equation}
holds locally informally in $x$.
A measurable and positive function $g: \mathbb{R}\to \mathbb{R}_+$ is said to be {\em self-neglecting} if \eqref{Rapid-eq-5} and \eqref{Rapid-eq-6} hold. 
 For example, if $g\in \mbox{RV}_\rho$ with $g(t)/t\to 0$, then the Uniform Convergence Theorem \ref{RV-the-1} implies that 
	\[\frac{g(t+xg(t))}{g(t)} = \frac{g\big(t(1+xg(t)/t)\big)}{g(t)}\to 1,
	\]  
	holds locally uniformly in $x$, thereby leading to that $g$ is self-neglecting. A less straightforward result is as follows.  
\begin{The} (Bloom \cite{Bloom76})
		\label{Rapid-the-0}
If $g$ is continuous and Beurling varying satisfying \eqref{Rapid-eq-5}, than $g$ is self-neglecting. 
\end{The}	

\begin{Cor}\label{Rapid-cor-0}
If a differentiable function $g:\mathbb{R}\to \mathbb{R}_+$ is ultimately monotone, satisfying that $g'(t)\to 0$, then $g$ is self-neglecting.
\end{Cor}

\noindent
{\sl Proof.} 
If $g$ is bounded from above, then $g(t)\to c$, where $c$ is a finite constant. This implies Beurling variation \eqref{Rapid-eq-5} and it follows from Theorem \ref{Rapid-the-0} that $g$ is self-neglecting. 

If $g$ is unbounded from above, then $g(t)\to \infty$. It follows from L'H\^{o}pital's rule that $g(t)/t\to 0$. Moreover, by the mean value theorem, since 
$g$ is differentiable, there exists a point $t^*\in [t,t+xg(t)]$ such that
\begin{equation}
	\label{Rapid-eq-61}
g(t+xg(t))-g(t) = g'(t^*)xg(t).
\end{equation}
Thus
\[\frac{g(t+xg(t))}{g(t)} = 1+xg'(t^*)\to 1,\ \mbox{as}\ t\to \infty. 
\]
By Theorem \ref{Rapid-the-0}, $g$ is self-neglecting. 
\hfill $\Box$

As a result, Beurling variation and self-neglecting properties are equivalent under the assumption of continuity. Since local uniform convergence plays a key role in representing rapid variation, the condition that 
$g$ is self-neglecting is often assumed in the rest of this paper.

\begin{The} (Omey \cite{Omey2013})
		\label{Rapid-the-1}
If $f\in \Gamma_\alpha(g)$ for some self-neglecting function $g$, then the convergence
\[\frac{f(t+xg(t))}{f(t)} \to \exp\{\alpha x\},\ \alpha\in \mathbb{R},
\]
holds locally informally in $x\in \mathbb{R}$. 
\end{The}

\begin{Rem}\label{Rapid-rem-2}
 The auxiliary function $g$ in $\Gamma_{\alpha}(g)$ is not unique, and the choice of $g$ can affect the rate of convergence of \eqref{Rapid-eq-4}. If $f\in \Gamma_{\alpha}(g)$, then  $f\in \Gamma_{\alpha}(g_1)$, where $g_1\sim g$, due to the uniform convergence described in Theorem \ref{Rapid-the-1}. 
\end{Rem}

If $f\in \Gamma_\alpha(g)$, where $g=1$, then it follows from  Theorem \ref{Rapid-the-1} and continuity of $\log t$, $h(t) := f(\log t)$ satisfies that $h(tx)/h(t)$ converges locally informally in $x$. This is true also due to Karamata's Theorem \ref{RV-the-1} (2), since $f(\log t)\in \mbox{RV}_\alpha$. 

On the other hand, if $f\in \mbox{RV}_\alpha$, then the Uniform Convergence Theorem \ref{RV-the-1} (2) implies that for any self-neglecting function $g$ such that $g(t)/t\to 0$, 
\[\frac{f(t+xg(t))}{f(t)} = \frac{f\big(t(1+xg(t)/t)\big)}{f(t)}\to 1,
\]
thereby resulting in $f\in \Gamma_0(g)$. That is, $\mbox{RV}_\alpha\subset \Gamma_0(g)$ with $g(t)/t\to 0$, and in particular, any regularly varying function $g$ with  $g(t)/t\to 0$ is self-neglecting and belongs to $\Gamma_0(g)$. 
Rapid variation is described by functions in $\Gamma_\alpha(g)$, $\alpha\ne 0$, and more precisely, if $\alpha>0$, then $\Gamma_\alpha(g)\subset \mbox{R}_\infty$, whereas if $\alpha<0$,  $\Gamma_\alpha(g)\subset \mbox{R}_{-\infty}$ \cite{de74}.

The auxiliary function $g$ in $\Gamma_{\alpha}(g)$ can be taken as the reciprocal of a hazard rate.

\begin{Pro}\label{Rapid-pro-1} Let $f:\mathbb{R}\to \mathbb{R}_+$ be differentiable with derivative $f'(t)$. 
	\begin{enumerate}
		\item If $f'(t)$ is ultimately positive and  $g(t)= f(t)/f'(t)$ is self-neglecting, then $f\in \Gamma_{1}(g)$. 
		\item If $f'(t)$ is ultimately negative and $g(t)= -f(t)/f'(t)$ is self-neglecting, then $f\in \Gamma_{-1}(g)$.
	\end{enumerate}
\end{Pro}

\noindent
{\sl Proof.} 
The proof can be found in \cite{Omey2013}. The proof is presented here to demonstrate that the convergence \eqref{Rapid-eq-4} depends on the behaviors of the underlying self-neglecting function $g$. 

(1) Consider, for $x>0$, 
\begin{equation}\label{Rapid-eq-71}
\log \frac{f(t+xg(t))}{f(t)}= \int_{t}^{t+xg(t)}\frac{f'(s)}{f(s)}ds=\int_{t}^{t+xg(t)}\frac{1}{g(s)}ds=\int_0^x\frac{g(t)}{g(t+zg(t))}dz.
\end{equation}
Since $g$ is self-neglecting, $\frac{g(t)}{g(t+zg(t))}\to 1$ uniformly on $[0,x]$, and it then follows that $\frac{f(t+xg(t))}{f(t)}\to e^x$.
The proof for $x<0$ is similar, and hence $f\in \Gamma_{1}(g)$.

(2)  Observe that 
\[\frac{d}{dt}\Big(\frac{1}{f(t)}\Big) = \frac{-f'(t)}{f^2(t)}
\]
is ultimately positive, and thus it follows from Proposition \ref{Rapid-pro-1} (1) above that since 
\[g(t) =\frac{-f(t)}{f'(t)}= \frac{1/f(t)}{-f'(t)/f^2(t)},
\]
is self-neglecting,  $1/f(t)\in \Gamma_{1}(g)$,  which implies that $f(t)\in \Gamma_{-1}(g)$.
\hfill $\Box$

\begin{Rem}\label{Rapid-rem-31}
	\begin{enumerate}
		\item It follows from Remark \ref{Rapid-rem-2} that under the assumption that $f'(t)$ is ultimately of one sign, if $g(t)f'(t)/f(t)\to \alpha$ and $g(t)$ is self-neglecting, then $f\in \Gamma_\alpha(g)$, for all $\alpha\ne 0$. 
		\item Note that \eqref{Rapid-eq-71} shows that the convergence \eqref{Rapid-eq-4} stems from the local uniform convergence of the underlying self-neglecting function $g(\cdot)$.
		\item   Suppose that $\overline{F}(t)$ denotes the survival function of a random variable with density $f(t)$ and $h(t) = f(t)/\overline{F}(t)$ denotes its hazard rate. If the reciprocal hazard rate $g(t) = 1/h(t)$ is self-neglecting, then $\overline{F}(t)\in \Gamma_{-1}(g)$. 
	\end{enumerate}
\end{Rem}

\begin{Exa}\rm\label{Rapid-exa-1}
\begin{enumerate}
	\item If $f(t) = \lambda e^{-\lambda x}$, $\lambda>0$, $x>0$, is the density of an exponential random variable, then $g(t)=1/\lambda$ is self-neglecting and thus $f\in \Gamma_{-1}(1/\lambda)$. 
	\item If $\varphi(t) = (2\pi)^{-1/2} e^{-t^2/2}$, $t\in \mathbb{R}$, is the density of a standard normal random variable, then $g(t)=-\varphi(t)/\varphi'(t) \sim 1/t$ is self-neglecting, and thus  $\varphi\in \Gamma_{-1}(1/t)$. 
	\item If $\overline{\Phi}(t) = \int_t^\infty (2\pi)^{-1/2} e^{-x^2/2}dx$, $t\in \mathbb{R}$, is the survival function of a standard normal random variable with the asymptotic
	hazard rate $h(t)\sim t$, then $\overline{\Phi}(t)\in \Gamma_{-1}(1/t)$. 
	\hfill $\Box$
\end{enumerate}
\end{Exa}

The auxiliary function $g$ in $\Gamma_{\alpha}(g)$ is taken as the reciprocal of a hazard rate in the rest of this paper. Define the cumulative hazard
\begin{equation}\label{Rapid-eq-7}
H(t) = \int_{t_0}^t \frac{1}{g(z)}dz = \int_{t_0}^t \frac{|f'(z)|}{f(z)}dz, \ t>t_0,
\end{equation}
for a differentiable function $f$ with derivative $f'$. Since $t/g(t)\to \infty$, $H(t)\to \infty$, as $t\to \infty$. Consider the total hazard construction $\Psi(t) = g(H^{-1}(t))$, which can be rewritten as
\[\Psi(t+x) = g\left(H^{-1}(t)+\frac{H^{-1}(t+x)-H^{-1}(t)}{g(H^{-1}(t))}g(H^{-1}(t))\right).
\]
If $g$ is self-neglecting, then by Theorem \ref{Rapid-the-1} and the uniform convergence, $\Psi(t+x)/\Psi(t)\to 1$, as $t\to \infty$, and thus $\Psi(\log t)\in \mbox{RV}_0$. It follows from Karamata's representation (see Theorem \ref{RV-the-2}) and simple substitutions that 
\begin{equation}\label{Rapid-eq-85}
	g(t) = C(t)\,\exp\Big\{\int_{t_0}^t\frac{\epsilon^*(z)}{g(z)}dz\Big\}
\end{equation}
where $C(t)\to c>0$, and $\epsilon^*(t) = \epsilon(e^{H(t)})$, where $\epsilon(t)\to 0$. Therefore, 
\begin{equation}\label{Rapid-eq-8}
g(t)\sim W(t) := c\,\exp\Big\{\int_{t_0}^t\frac{\epsilon^*(z)}{g(z)}dz\Big\}
\end{equation}
where $\epsilon^*(t)\to 0$, and $c>0$ is a constant, leading to $W'(t)\sim \epsilon^*(t)\to 0$. It turns out that \eqref{Rapid-eq-8} is also a sufficient condition for $g(\cdot)$ to be self-neglecting, as shown in the following result from \cite{Omey2013}.

\begin{The}	\label{Rapid-the-2}
A measurable and positive function $g$ is self-neglecting if and only if $g(t)\sim W(t)$ with $W'(t)\to 0$. 
\end{The}

\begin{Rem}\label{Rapid-rem-3}
\begin{enumerate}
	\item Theorem \ref{Rapid-the-2} strengthens Corollary \ref{Rapid-cor-0}, and hereafter a self-neglecting function $g$ can be taken as  the reciprocal of a hazard rate whose derivative converges to zero. 
	\item A measurable and positive function $g(t)$ is self-neglecting if and only if $g(H^{-1}(\log t))$ is slowly varying, where $H(\cdot)$ is the cumulative hazard \eqref{Rapid-eq-7}. 
\end{enumerate}
\end{Rem}

\begin{The}\label{Rapid-the-3}
	(Omey \cite{Omey2013}) 
Let $f: \mathbb{R}\to \mathbb{R}_+$ be measurable. A function $f\in \Gamma_\alpha(g)$, where $g$ is self-neglecting, if and only if 
\begin{equation}\label{Rapid-eq-9}
f(t) \sim A(t) \exp\Big\{\alpha\int_{t_0}^t\frac{1}{g(z)}dz\Big\}, \ A(t)=\exp\{C+w(t)\}, 
\end{equation}
for a constant $C$ and  $w(t) = \int_{t_0}^t\frac{\epsilon(z)}{g(z)}dz$, where $\epsilon(\cdot)$ is given by \eqref{Rapid-eq-85}, such that $A(\cdot)$  satisfies that $g(t)A'(t)/A(t)\to 0$. 
\end{The}

That is, any positive function $f\in \Gamma_\alpha(g)$ can be rephrased as $f(t) \sim A(t) h(t)$, where $h(t) = \exp\big\{\alpha\int_{t_0}^tdz/g(z)\big\}$, satisfying that 
\[g(t)A'(t)/A(t)\to 0,\ g(t)h'(t)/h(t)\to \alpha, \ \alpha\in \mathbb{R}. 
\]
Using Bachmann–Landau notations, $|A'(t)|/A(t) = o(1/g(t))$. Therefore, the hazard rate of a function $A(\cdot)$ is dominated by the hazard rate $1/g(t)$ asymptotically if and only if $A(t)\in \Gamma_0(g)$. The fact that $\mbox{RV}_\alpha\subset \Gamma_0(g)$ can be easily explained by Karamata's theorem, since the hazard rate of any regularly varying function is dominated asymptotically by $1/g(t)$, where $g(t)$ is self-neglecting.

Moreover, combining two functions $w(\cdot)$ and $A(\cdot)$, \eqref{Rapid-eq-9} yields the following representation \cite{de74}.

\begin{The}\label{Rapid-the-37}
	Let $f: \mathbb{R}\to \mathbb{R}_+$ be measurable. A function $f\in \Gamma_\alpha(g)$, where $g$ is self-neglecting, if and only if 
	\begin{equation}\label{Rapid-eq-99}
		f(t) \sim B(t)h(t),\ \mbox{where}\  h(t)=  \exp\Big\{\int_{t_0}^t\frac{\alpha+\epsilon(z)}{g(z)}dz\Big\}, 
	\end{equation}
	where $B(t)>0$ with $\lim_{t\to \infty}B(t)=b>0$ and $\lim_{t\to \infty}\epsilon(t)=0$, such that $g(t)h'(t)/h(t)\to \alpha$. 
\end{The}

\begin{Rem}\label{Rapid-rem-35}
	\begin{enumerate}
		\item The representation \eqref{Rapid-eq-99} can be rewritten as 
		\begin{equation}\label{Rapid-eq-155}
			f(t)\sim b\exp\Big\{\int_{t_0}^t\frac{(\alpha+\epsilon(z))\lambda(z)}{z}dz\Big\},\ \alpha\in \mathbb{R}, 
		\end{equation}
	where $b>0$ and $\lambda(t)\sim t/g(t)$, as $t\to \infty$. That is, $\lambda(\cdot)$ goes to infinity at the same rate as that of $t/g(t)$, which is in contrast to Karamata's representation for regular variation (Theorem \ref{RV-the-2}). 
		\item For any regularly varying function $h(t)\in \mbox{RV}_\rho$, Karamata's conditions \eqref{RV-eq-7} and \eqref{RV-eq-8} can be rephrased as $t\,h'(t)/h(t)\to \rho$, $\rho\in \mathbb{R}$.
		In contrast,  any $f\in \Gamma_\alpha(g)$ can be rephrased as $f(t) \sim B(t) h(t)$, satisfying that 
		\[g(t)h'(t)/h(t)\to \alpha, \ \alpha\in \mathbb{R}, \mbox{and}\ t\,h'(t)/h(t)\to \infty. 
		\]
		\item In general, asymptotic tail variation can be written, in terms of hazard rates, as follows. 
		\begin{equation}
			\lim_{t\to \infty}t\,\frac{h'(t)}{h(t)}=	\left\{\begin{array}{ll}
				0 & \mbox{(slow variation at $\infty$)}\\
				\rho\in \mathbb{R}\backslash \{0\}  & \mbox{(regular variation at $\infty$)}\\
				\pm \infty & \mbox{(rapid variation at $\infty$).}
			\end{array}
			\right.
			\label{Rapid-eq-15}
		\end{equation}
	Note that $\mbox{RV}_{\rho}\subset \Gamma_0(g)$ with $g(t)/t\to 0$, and $\Gamma_\alpha(g)\subset R_{\infty}\cup R_{-\infty}$, $\alpha\ne 0$.  
	\end{enumerate}
\end{Rem}

Utilizing regularly varying functions $g$ with $g(t)/t\to 0$ in \eqref{Rapid-eq-9}, typical functions in $\Gamma_{\alpha}(g)$ have forms that involve terms such as $r(t)\exp\big\{\alpha t^b\big\}$, where $\alpha\in \mathbb{R}$, $b>0$ and $r(\cdot)$ is regularly varying.

\begin{Cor}\label{Rapid-cor-4}
If $f\in \Gamma_{\alpha}(g)$, where $g$ is ultimately decreasing and $\alpha<0$ (increasing and $\alpha>0$), then $f$ is tail equivalent to a function that is ultimately decreasing (increasing).  
\end{Cor}

\noindent
{\sl Proof.} If $g$ is ultimately decreasing, and $\alpha<0$, then it follows from \eqref{Rapid-eq-85} and \eqref{Rapid-eq-8} that $w(t)$ is ultimately decreasing, and therefore, by Theorem \ref{Rapid-the-3}, $f$ is tail equivalent to a function that is ultimately decreasing. The proof of the other case is similar. 
\hfill $\Box$

Consider now a random variable  $X$ with survival function $\overline{F}(x)$, density $f(x)$ and hazard rate $h(t) = f(t)/\overline{F}(x)$. Let $g(t) = \overline{F}(x)/f(t)$ denote the reciprocal of the hazard rate of $X$. According to the von Mises condition (see \eqref{RV-eq-20}), if $t\,h(t) = t/g(t)$ converges to a finite limit, then $X$ is regularly varying at $\infty$. In contrast, rapid variation emerges from the situation that  $t\,h(t) = t/g(t)$ converges to infinity. See also Remark \ref{Rapid-rem-35}.

\begin{Cor}\label{Rapid-cor-35}
Let $X$ be a random variable with survival function $\overline{F}(x)$ and density $f(x)$. The reciprocal of the hazard rate $1/h(x) = \overline{F}(x)/f(x)$ is self-neglecting, if and only if, $\overline{F}(x)\in \Gamma_{-1}(h^{-1})$.
\end{Cor}

\noindent
{\sl Proof.} It follows from Theorem \ref{Rapid-the-3} and the fact that
$\overline{F}(x)=\exp\{-\int_{-\infty}^xh(t)dt\}$, $x\in \mathbb{R}$. 
\hfill $\Box$

It is worth mentioning that 
as discussed for regular variation, $\mbox{RV}_\alpha\subseteq \Gamma_0(g)$, and rapid variation is described by the class $\Gamma_\alpha(g)$, $\alpha\ne 0$. If a survival function $\overline{F}(x)$ is regularly varying, then $\overline{F}(x)\in \Gamma_0(g)$, whereas for $\overline{F}(x)\in \Gamma_{-1}(g)$, where $g(t) = 1/h(t)$ is the reciprocal of the hazard rate, 
\[\lim_{t\to \infty}\frac{\overline{F}(tx)}{\overline{F}(t)} = \exp\Big\{-\int_{t}^{tx}h(z)dz\Big\}, 
\]
which shows that $\overline{F}(x)$ is rapidly varying in the sense of de Haan \eqref{Rapid-eq-1}.

\begin{Exa}\rm\label{Rapid-exa-2}
Let $X$ be a positive random variable with hazard rate $h(t) = t^k$, $k>-1$, $t> 0$. The reciprocal of the hazard rate $g(t) = t^{-k}$, and $g'(t) = -kt^{-k-1}\to 0$, as $t\to \infty$. It follows from Corollary \ref{Rapid-cor-35} that the survival function $\overline{F}(x)\in \Gamma_{-1}(t^{-k})$, $k>-1$. Note that if $k\ge 0$, $X$ has the increasing hazard rate, and if $0>k>-1$, then $X$ has the decreasing hazard rate. 
\hfill $\Box$
\end{Exa}

 Some useful properties are summarized below.

\begin{Pro}
	\label{Rapid-pro-2}
	\begin{enumerate}
		\item If $f\in \Gamma_\alpha(g)$, $\alpha\in \mathbb{R}$, then $f^\beta\in \Gamma_{\alpha\beta}(g)$, $\forall\ \beta\in \mathbb{R}$. 
		\item If $f\in \Gamma_\alpha(g)$, $\alpha\in \mathbb{R}$, then for any $\epsilon>0$ and a compact subset $B\subset \mathbb{R}$, there exists $t_0=t_0(\epsilon, B)$, such that whenever $t\ge t_0$, 
		\begin{equation}
			\label{Rapid-eq-55}
			(1-\epsilon)e^{\alpha x}<\frac{f(t+xg(t))}{f(t)}< (1+\epsilon)e^{\alpha x},
		\end{equation}
	for all $x\in B$. 
		\item If $f_1\in \Gamma_{\alpha_1}(g)$ and $f_2\in \Gamma_{\alpha_2}(g)$, $\alpha_1, \alpha_2\in \mathbb{R}$, than $f_1f_2\in \Gamma_{\alpha_1+\alpha_2}(g)$. 
	\end{enumerate}
\end{Pro}

Similar to Karamata's thereon (Theorem \ref{RV-the-3}), the following relation between a function and its derivative with respect to the $\Gamma_\alpha(g)$ class can be established. 

\begin{The}\label{Rapid-the-4}
	\begin{enumerate}
		\item Suppose that $f$ has an ultimately increasing (decreasing) derivative $f'$. If $f\in \Gamma_\alpha(g)$, then $|f'|\in \Gamma_\alpha(g)$. Moreover, 
		\begin{enumerate}
			\item If $f\in \Gamma_{\alpha}(g)$, $\alpha>0$, has an ultimately increasing (decreasing) derivative $f'$, then $g(t)\sim f(t)/f'(t)$ is self-neglecting. 
			\item If $f\in \Gamma_{\alpha}(g)$, $\alpha<0$, has an ultimately increasing (decreasing) derivative $f'$, then $g(t)\sim -f(t)/f'(t)$ is self-neglecting. 
		\end{enumerate}
		\item Suppose that $f$ is ultimately monotone. If $f$ has a derivative $f'$, such that $|f'|\in \Gamma_\alpha(g)$, then $f\in \Gamma_\alpha(g)$.
	\end{enumerate}
\end{The}

\noindent
{\sl Proof.} (1) The proof of the case that $f'$ is ultimately increasing can be found in \cite{Omey2013}, and the similar proof for $\Gamma_{\alpha}(g)$, $\alpha<0$, is outlined here to illustrate the von Mises-type condition for rapid variation. Without loss of generality, assume that $\alpha = -1$. 

Observe that 
\begin{equation}\label{Rapid-eq-16}
f(t+xg(t))-f(t)=g(t)\int_0^xf'(t+zg(t))dz,\ x>0.
\end{equation}
For sufficiently large $t$, since $f'(t)$ is increasing, then
\[xg(t)f'(t)\le f(t+xg(t))-f(t)\le xg(t)f'(t+xg(t))
\]
leading to 
\[
\frac{xg(t)f'(t)}{f(t)}\le \frac{f(t+xg(t))}{f(t)}-1\le \frac{xg(t)f'(t+xg(t))}{f(t)}.
\]
Taking the limit on the first inequality, one has
\begin{equation}\label{Rapid-eq-14}
\limsup_{t\to \infty}\frac{g(t)f'(t)}{f(t)}\le  \frac{1}{x}(e^{- x}-1).
\end{equation}
Rewrite the second inequality as 
\begin{equation}\label{Rapid-eq-13}
\frac{g(t+xg(t))}{g(t)}\frac{f(t)}{f(t+xg(t))}\left(\frac{f(t+xg(t))}{f(t)}-1\right)\le \frac{xg(t+xg(t))f'(t+xg(t))}{f(t+xg(t))}.
\end{equation}
Since $f\in \Gamma_{\alpha}(g)$, where $g$ is self-neglecting, then $g(t)\sim g(t+xg(t))=g(t_1)$. Taking the limit on \eqref{Rapid-eq-13} leads to
\begin{equation}\label{Rapid-eq-15}
\frac{1}{x}e^{ x}(e^{- x}-1)\le \liminf_{t\to \infty}\frac{g(t)f'(t)}{f(t)}.
\end{equation}
As $x\to 0$, it follows from \eqref{Rapid-eq-14} and \eqref{Rapid-eq-15} that 
$\lim_{t\to \infty}\frac{g(t)f'(t)}{f(t)}=-1$. Since $f\in \Gamma_{\alpha}(g)$, where $g(t) \sim -f(t)/f'(t)$ is self-neglecting, it then  follows from the representations \eqref{Rapid-eq-9} and \eqref{Rapid-eq-85} that $-f'\in \Gamma_{\alpha}(g)$.

(2) Consider, without loss of generality, that $f$ is ultimately decreasing and  $-f'\in \Gamma_{-1}(g)$. The proof of the other case that $f'\in \Gamma_{1}(g)$ is similar.

Since $-f'\in \Gamma_{-1}(g)$, it follows from Proposition \ref{Rapid-pro-2} (2) that for any $\epsilon>0$, there exists $t_0(\epsilon, x)$ such that whenever $t>t_0(\epsilon, x)$, 
\[	(1-\epsilon)e^{- z}<\frac{f'(t+zg(t))}{f'(t)}< (1+\epsilon)e^{- z},\ \forall\ z\in [0,x].
\]
Integrating with respect to $z$ yields
\[(1-\epsilon)(1-e^{-x})g(t)<\frac{f(t+xg(t))-f(t)}{f'(t)}<(1+\epsilon)(1-e^{-x})g(t),
\]
leading to
\begin{equation}\label{Rapid-eq-155}
(1-\epsilon)(1-e^{-x})\frac{g(t)f'(t)}{f(t)}>\frac{f(t+xg(t))-f(t)}{f(t)}>(1+\epsilon)(1-e^{-x})\frac{g(t)f'(t)}{f(t)}.
\end{equation}
Take the limit as $t\to \infty$ on the first inequality, one has the following
\[(1-\epsilon)(1-e^{-x})\liminf_{t\to \infty}\frac{g(t)f'(t)}{f(t)}\ge \liminf_{t\to \infty}\frac{f(t+xg(t))-f(t)}{f(t)}.
\]
Let $\epsilon\to 0$ and then $x\to \infty$, and it follows from the ultimate monotonicity and Dini's theorem that 
\[\liminf_{t\to \infty}\frac{g(t)f'(t)}{f(t)}\ge -1.
\]
Similarly, take the limit as $t\to \infty$ on the second inequality of \eqref{Rapid-eq-155}, one also has
\[\limsup_{t\to \infty}\frac{f(t+xg(t))-f(t)}{f(t)}\ge (1+\epsilon)(1-e^{-x})\limsup_{t\to \infty}\frac{g(t)f'(t)}{f(t)}.
\]
Let $\epsilon\to 0$ and then $x\to \infty$, and it follows again from the ultimate monotonicity and Dini's theorem that 
\[-1 \ge \limsup_{t\to \infty}\frac{g(t)f'(t)}{f(t)}.
\]
Therefore, $\lim_{t\to \infty}\frac{-g(t)f'(t)}{f(t)}=1$. Since $g(t)\sim -f(t)/f'(t)$ is self-neglecting, by Proposition \ref{Rapid-pro-1} (2), $f\in \Gamma_{-1}(g)$. 
\hfill $\Box$

A Taylor expansion for $f$ at $\infty$ can be also established. 
\begin{Pro}\label{Rapid-pro-3}
If $f$ has a derivative $f'\in \Gamma_0(g)$, then the first-order Taylor expansion holds at $\infty$:
\[f(t+xg(t))-f(t) \sim xg(t)f'(t), \ x\in \mathbb{R}. 
\]
\end{Pro}

\noindent
{\sl Proof.}
It follows from \eqref{Rapid-eq-16} that 
\[\frac{f(t+xg(t))-f(t)}{g(t)f'(t)}=\int_0^x\frac{f'(t+zg(t))}{f'(t)}dz
\]
Since $f'\in \Gamma_0(g)$, then, by Theorem \ref{Rapid-the-1}, $\frac{f(t+xg(t))-f(t)}{g(t)f'(t)}\to x$ and the result follows. 
\hfill $\Box$

\begin{Rem}\label{Rapid-rem-31}
	\begin{enumerate}
			\item If $f\in \Gamma_{\alpha}(g)$ has an ultimately monotone derivative $f'$, then the von Mises-type condition for rapid variation can be written as $\lim_{t\to \infty}\frac{g(t)f'(t)}{f(t)}=\alpha$, $\alpha\in \mathbb{R}$. Specifically, if $f:\mathbb{R}\to \mathbb{R}$ is ultimately positive with ultimately monotone derivative $f'(t)$, then
			\begin{equation}\label{Rapid-eq-29}
			\lim_{t\to \infty}\frac{g(t)f'(t)}{f(t)}=\alpha\in \mathbb{R}\ \mbox{and $g(\cdot)$ is self-neglecting}\ \mbox{if and only if}\ f\in \Gamma_\alpha(g).
		\end{equation}
			 See also Theorem \ref{Rapid-the-3} and Remark \ref{Rapid-rem-35} (2) for the von Mises-type condition, in terms of the tail equivalence, for rapid variation.
		\item If a differentiable function $g:\mathbb{R}\to \mathbb{R}_+$ is ultimately monotone, satisfying that $g'(t)\to 0$, then, by Corollary \ref{Rapid-cor-0}, $g$ is self-neglecting. Therefore, $|g'| \in \Gamma_0(g)$, and 
		\[g(t+xg(t))-g(t) \sim xg(t)g'(t),
		\]
	which strengthens \eqref{Rapid-eq-61}.  
	\end{enumerate}
\end{Rem}	

If $f$ is regularly varying with tail parameter $\rho$, $\rho\ne 0$, then the inverse of $f$ is regularly varying with tail parameter $1/\rho$. The inverse function should be understood as left or right continuous inverse that is consistent with the ultimate monotonicity. The inverses of functions from $\Gamma_\alpha(g)$, $\alpha\ne 0$, do not belong to the same class 
 but form the class $\Pi$, which is a sub-class of slow variation, introduced and studied by de Haan \cite{de74}. Consider first $\Gamma_{\alpha}(g)$, $\alpha>0$, consisting of the functions that are ultimately increasing.

 \begin{Def}\label{Rapid-def-2}
 An ultimately increasing function $v: \mathbb{R}_+\to \mathbb{R}$, with $\lim_{x\to \infty}v(x)=\infty$, is said to be in the class $\Pi(a,b)$, denoted by $v\in \Pi(a,b)$, if there exist functions $a: \mathbb{R}_+\to \mathbb{R}_+$ and $b:\mathbb{R}_+\to \mathbb{R}$, such that for any $x>0$, 
 \begin{equation}\label{Rapid-eq-20}
\lim_{t\to \infty}\frac{v(tx)-b(t)}{a(t)}= \log x.
 \end{equation}
 \end{Def}

\begin{Rem}\label{Rapid-rem-5}
\begin{enumerate}
	\item The class of auxiliary functions $a(t)$ is $\mbox{RV}_0$. 
	\item The convergence \eqref{Rapid-eq-20} holds locally uniformly in $x$. 
	\item The representation of $v\in \Pi(a,b)$ is given by
	\begin{equation}\label{Rapid-eq-21}
	v(x) = h(x) + \int_{x_0}^x h(t)dt/t, 
	\end{equation}
for a slowly varying function $h$. 
\end{enumerate}
\end{Rem}

If $h\in \mbox{RV}_0$, then $t^{-1}h(t)\in \mbox{RV}_{-1}$. Theorem \ref{RV-the-3} (1) implies that $\int_{x_0}^x h(t)dt/t\in \mbox{RV}_0$. Therefore, it follows from \eqref{Rapid-eq-21} that $\Pi(a,b)\subset \mbox{RV}_0$. 

\begin{Pro}\label{Rapid-the-5} (de Haan \cite{de74}) 
Suppose that  $f: \mathbb{R} \to \mathbb{R}_+$ is an ultimately increasing function, with $\lim_{t\to \infty}f(t)=\infty$. Then $f\in \Gamma_{\alpha}(g)$, with $\alpha>0$, if and only if the left-continuous inverse $f^{\leftarrow}\in \Pi(a,b)$, where $a(t) = g(f^{\leftarrow}(t))$.
\end{Pro}

\begin{Rem}\label{Rapid-rem-6}
	\begin{enumerate}
		\item For $f\in \Gamma_{\alpha}(g)$, with $\alpha>0$, $f^{\leftarrow}(x) = a(x)+\int_{x_0}^xa(t)dt/t\in \Pi(a,b)$, where one can choose that 
		\[a(t) = f(t)-t^{-1}\int_{t_0}^tf(z)dz
		\]
		and $b(t)=f(t)$. 
		\item In general, $f\in \mbox{R}_\infty$ if and only if $f^{\leftarrow}\in \mbox{RV}_0$ \cite{Elez2013}. 
	\end{enumerate}
\end{Rem}

Similar to Proposition \ref{RV-pro-25} (6), the tail equivalence related to rapid variation is described as follows.

\begin{Pro}\label{Rapid-the-6}
	 The two functions $f_1, f_2\in \Gamma_{\alpha}(g)$, with $\alpha>0$, if and only if the left-continuous inverses $f_1^{\leftarrow}(x) = f_2^{\leftarrow}(U(x))$, for some $U\in \mbox{RV}_1$. 
\end{Pro}

Suppose that $f: \mathbb{R} \to \mathbb{R}_+$ is an ultimately decreasing function, with $\lim_{t\to \infty}f(t)=0$. If$f\in \Gamma_{\alpha}(g)$, with $\alpha<0$, then $\hat{f}(t)=1/f(t)\in \Gamma_{-\alpha}(g)$ is ultimately increasing, with $\lim_{t\to \infty}\hat{f}(t)=\infty$. Observe that 
\[\hat{f}^{\leftarrow}(t^{-1}) = \inf\big\{s: \hat{f}(s)\ge t^{-1}\big\}=\inf\big\{s: f(s)\le t\big\}=f^{\rightarrow}(t),
\]
and then the right-continuous inverse $f^{\rightarrow}(t^{-1})\in \Pi(a,b)$, where
	\begin{equation}\label{Rapid-eq-22}
	a(t) = \frac{1}{f(t)}-t^{-1}\int_{t_0}^t\frac{1}{f(z)}dz\in \mbox{RV}_0,\ b(t) = 1/f(t)\in \Gamma_{-\alpha}(g). 
\end{equation}

\begin{Pro}\label{Rapid-the-6}
	Suppose that  $f: \mathbb{R} \to \mathbb{R}_+$ is an ultimately decreasing function, with $\lim_{t\to \infty}f(t)=0$. Then $f\in \Gamma_{\alpha}(g)$, with $\alpha<0$, if and only if the right-continuous inverse $f^{\rightarrow}(t^{-1})\in \Pi(a,b)$, or equivalently for any $x>0$, 
	\begin{equation}\label{Rapid-eq-24}
		\lim_{u\to 0}\frac{f^{\rightarrow}(ux)-b(u^{-1})}{a(u^{-1})}= -\log x,
	\end{equation}
for some functions $a(\cdot)$ and $b(\cdot)$, satisfying \eqref{Rapid-eq-22}. 
\end{Pro}

\begin{Rem}\label{Rapid-rem-7}
\begin{enumerate}
	\item The convergence \eqref{Rapid-eq-24} holds locally uniformly in $x$.
	\item The right-continuous inverse $f^{\rightarrow}(u)$ is known as being  slowly varying at $0$, and can be represented as
		\begin{equation}\label{Rapid-eq-23}
		f^{\rightarrow}(u) = a(u^{-1}) - \int_{u_0}^u a(z^{-1})dz/z,  
	\end{equation}
where $a(\cdot)$ is given in \eqref{Rapid-eq-22}. 
\item In general, $f\in \mbox{R}_{-\infty}$ if and only if $f^{\rightarrow}(u)$ is slowly varying at $0$ \cite{Elez2013}. 
\end{enumerate}
\end{Rem}

\begin{Pro}\label{Rapid-the-9}
	The two functions $f_1, f_2\in \Gamma_{\alpha}(g)$, with $\alpha<0$, if and only if the right-continuous inverses $f_1^{\rightarrow}(x) = f_2^{\rightarrow}(V(x))$, for some $V(\cdot)$ is regularly varying at 0 with tail parameter $1$. 
\end{Pro}

It is worth mentioning that an extended Gamma class of rapidly varying functions is introduced and studied in \cite{Omey2013}, and this class contains both $\Gamma_{\alpha}(g)$ and $\Pi(a,b)$ as two sub-classes.

\section{Concluding Remarks}
\label{univeriate conclusion}

The representation for functions with various tail variations are summarized as follows.
\begin{enumerate}
	\item A function $f\in \mbox{RV}_\alpha$ if and only if $f(t)\sim 
	\exp\left\{c+\int_{t_0}^t\frac{A(z)}{z}dz\right\}$, where $A(t)\to \alpha$, as $t\to \infty$. 
	\item A function $f\in \Gamma_\alpha(g)$ if and only if $f(t)\sim 
	\exp\left\{c+\int_{t_0}^t\frac{A(z)}{g(z)}dz\right\}\sim 
	\exp\left\{c+\int_{t_0}^t\frac{A(z) B(z)}{z}dz\right\}$, where $A(t)\to \alpha$ and $B(t)\sim t/g(t)\to \infty$, as $t\to \infty$. 
	\item A function $f\in R_{\pm \infty}$ if and only if $f(t)\sim 
	\exp\left\{c+\int_{t_0}^{T(t)}\frac{A(z)}{z}dz\right\}$, where $T(t)\sim t$ and $A(t)\to \pm\infty$, as $t\to \infty$.
\end{enumerate}
The classes $R_{\pm \infty}$ of rapid variation in the sense of de Haan are studied extensively in \cite{Elez2013, Elez2015}, from which the following
characterizations for regular and rapid variation can also be made in terms of Karamata's forms.
\begin{enumerate}
	\item A function $f\in \mbox{RV}_\alpha$ if and only if $f(t) = t^\alpha \ell(t)$, where $\ell(tx)\sim \ell(t)$ for $x>0$.
	\item A function $f\in R_{\infty}$ ($R_{-\infty}$) if and only if $f(t) = t^\alpha A_\alpha(T_\alpha(t))$, for all $\alpha>0$ ($<0$), where $T_\alpha(t)\sim t$ and $A_\alpha: (0,\infty)\to (0,\infty)$ is increasing (decreasing). 
\end{enumerate} 

Note that  $g(t)$ is self-neglecting if and only if $\ell(t)=g(H^{-1}(\log t))$ is slowly varying, where $H(\cdot)$ is the cumulative hazard \eqref{Rapid-eq-7}.
Therefore, the classes $\mbox{RV}_\alpha$ and $\Gamma_\alpha(g)$ are well-defined, with  self-neglecting functions $g(\cdot)$, interpreted as reciprocal hazard rates,  play a central role in their analysis.  In contrast, the challenge with the classes $R_{\pm \infty}$ arises from the behavior of the function $A(\cdot)$, which diverges to either positive or negative infinity without a specified rate of convergence. Nevertheless, the $R_{\pm \infty}$ classes encompass a broader range of functions, including highly super-increasing and super-decreasing behaviors.

\end{document}